\documentclass[12pt]{article}

\usepackage{amsmath,amssymb}
\usepackage{verbatim}
\usepackage{comment}
\usepackage{color}
\usepackage{amsthm}
\usepackage{cases}
\usepackage{url}

\newtheorem{theorem}{Theorem}[section]
\newtheorem{thm}[theorem]{Theorem}

\newtheorem{cor}[theorem]{Corollary}
\newtheorem{lemma}[theorem]{Lemma}

\newtheorem{prop}[theorem]{Proposition}
\newtheorem{definition}[theorem]{Definition}
\newtheorem{defn}[theorem]{Definition}
 
\newtheorem{assumption}[theorem]{Assumption}
\newtheorem{remark}[theorem]{Remark}
\newtheorem{example}[theorem]{Example}
\numberwithin{equation}{section}

\def\be{\begin{equation}}
\def\ee{\end{equation}}
\def\bes{\begin{equation*}}
\def\ees{\end{equation*}}
\def\Mon{{\rm Mon}}
\def\Tue{{\rm Tue}}
\def\Heads{{\rm Heads}}
\def\Tails{{\rm Tails}}

\def\PA{{P_A}}

\def\FOb{{\sF_W}}
\def\Ber{{\rm Ber}}

  \def\sF {{\mathcal F}}
  
 \def\sK {{\mathcal K}} 
  
\def\sP {{\mathcal P}}  
  
  \def\sX {{\mathcal X}}
 \def\sZ {{\mathcal Z}}

 \def\bE {{\mathbb E}}

 \def\bN {{\mathbb N}} 
\def\bP {{\mathbb P}}  \def\bR {{\mathbb R}}

 \def\bZ {{\mathbb Z}}

\def\sms{\smallskip}
\def\ms{\medskip}

\def\sm{\smallskip\noindent}

\def\ignore#1{}

\def\al {\alpha}
\def\lam {\lambda} \def\Lam {\Lambda} 
\def\eps{\varepsilon}
\def\th{\theta} \def\Th{\Theta}

\def\Om{\Omega} \def\om{\omega}

\def\to {\rightarrow}

\def\pd {\partial}
\def\q{\quad} 
\def\dint{\int\kern-.6em\int}

\def \half {{\tfrac12}}
\def\fract{\tfrac}
\def \third {{\tfrac13}}       

\def\wt{\widetilde}

\def\be{\begin{equation}}
\def\ee{\end{equation}}
\def\bes{\begin{equation*}}
\def\ees{\end{equation*}}
\def\ba{\begin{align}}
\def\ea{\end{align}}
\def\xxea{\end{align}}
\def\bas{\begin{align*}}
\def\eas{\end{align*}}

\def\proof{{\smallskip\noindent {\em Proof. }}}
\def\qed{{\hfill $\square$ \bigskip}}

\definecolor{dgreen}{rgb}{0, 0, 0.1}
\definecolor{dblue}{rgb}{0, 0.0, 0.6}
\definecolor{vdblue}{rgb}{0,.08, 0.45}
\definecolor{dred}{rgb}{0.7, 0.0, 0.0}
\definecolor{vdblue}{rgb}{0,.08, 0.45}

\definecolor{purple}{rgb}{0.6, 0.0, 0.6}
\definecolor{mytext}{rgb}{0.1, 0.1, 0.1}

\begin{document}

\font\titlefont=cmbx14 scaled\magstep1
\title{\titlefont  \vspace{-5ex}  Inference by Multiple Identical Observers } 

\author{
M. T. Barlow\footnote{ Department of Mathematics,
University of British Columbia,
Vancouver, BC V6T 1Z2, Canada. barlow@math.ubc.ca
}
}

\maketitle


\begin{abstract}
We consider models for inference made by observers which may have multiple identical copies,
such as in the Sleeping Beauty problem.
We establish a framework for describing these problems on a probability space satisfying Kolmogorov's axioms.
This enables the main competing approaches to be compared precisely.
\end{abstract}

\section{Introduction} 

This paper is concerned with inference by observers whose existence and multiplicity is affected by the outcome
of the random experiment which they are observing. 
One can argue that implicit in the classical model of a random experiment is the assumption that 
observations and statistical inferences are made by a Classical  Observer (CO), 
which exists before the start of the experiment, and whose 
existence is not affected by the outcome of the experiment. 
Little attention has been paid to this requirement since these conditions are nearly
always satisfied in practice.

However non Classical Observers do arise in physics and cosmology. 
The laws of nature depend on about 20-30 
constants, many of which appear ``fine tuned' for life -- see
\cite{Ad19, MR}. Among these is the cosmological constant $\Lam$, which has a value
much smaller than that predicted by some theories. One type of explanation 
is via a multiverse, which contains many universes, each with a different value of $\Lam$.
Only a small range of these values will permit the existence of stars, planets, life and observers.
The value for our universe $\Lam_0$ has to lie in this range, and we can ask whether
this observed value is likely given the multiverse theory. To answer this we need to know
how the (random) observers in these universes should make inferences. 
(For some recent papers in this area see \cite{BES18, SPL}.)

Following the terminology introduced by Brandon Carter in \cite{BC}, I
will call the conditional observers arising in examples like this 
{\em Anthropic Observers} or AO. (The term is flawed as these observers do not have to be human.) 
Most examples allow multiple
observers, and one needs to know how  the multiplicity of these observers affects  inferences. 

A `toy model' one can use to  investigate these questions is the Sleeping Beauty problem,
described by A. Elga in \cite{E} in the philosophy journal {\em Analysis}.
The name is due to Robert Stalnaker,  an earlier version of the
problem is due to A. Zuboff (see \cite{Z}),  and an essentially identical problem is given in Example 5 of \cite{PR}, a paper on the Absentminded Driver problem in decision theory. 
Elga introduced the problem in order to examine 
how observers update probabilities  in centered and uncentered worlds.

On Sunday the experiment, as described below,
is explained to Sleeping Beauty (SB). On Sunday night she goes to sleep in an isolated cell. 
A fair coin is tossed. Whatever the outcome she is awakened on Monday morning, and then on Monday
evening is given a potion which puts her to sleep and causes her to forget everything that occurred
on Monday. 
If the outcome of the toss was Tails she is then woken on Tuesday, while if the outcome was Heads
then she  is not woken on Tuesday.  
We ask two questions:

\sm {\em Question 1.} 
 When SB wakes in her cell during the experiment, what probability should she assign
 to the event that the coin landed Heads? 

 \sm {\em Question 2.} 
 If she is then told that it is Monday, what probability should she now assign to Heads? 
 
 The Sleeping Beauty problem has been framed in a fanciful fashion, 
 but it raises questions which are relevant
in mainstream research in physics and cosmology. These problems have been actively
discussed by physicists and philosophers, but they also deserve attention from the
mathematical community.
 
 Philosophers' use the term {\em credence} to denote a subjective probability.
 I will use that term as  shorthand for ``the probability law that SB should use'', but
 without wishing to be tied to any particular philosophical interpretation of probability.
 There are four positions in the philosophy literature. {\em Halfers} claim that the answer to
 Question 1 is $\half$. 
 {\em Standard Halfers} claim  that the answer to Question 2 is $\fract23$, while 
 {\em Double Halfers} argue that the answer to Question 2 is still $\half$. 
 {\em Thirders} argue that the answer to Question 1 is $\fract13$, and then deduce that
 the answer to Question 2 is $\half$. A final group argue that the problem is ambiguous,
 inappropriate or undetermined. 
 In his excellent review \cite{W} Winkler suggests that the Thirder view is the majority position among philosophers, but that this is not reflected in the literature due to publication bias.
Physics papers such as \cite{BES18} use a weighting which aligns with the Thirder approach.

 A further question is whether this problem can be adequately described using 
Kolmogorov's axioms for probability. 
The main goal of this paper is to argue that it can. 
The  framework we construct then allows the competing views of Halfers and Thirders 
to be evaluated in precise mathematical terms. 

\ms
Elga  \cite{E} proposed a straightforward solution of the SB problem using conditional probabilities. 
He claimed that
\begin{align}
  \label{e:EPI}
  \bP( {\rm Mon}  | {\rm Tails} )  & = \half, \\ 
 \label{e:EPEI}
  \bP( {\rm Heads}  | {\rm Mon } )  &=\bP( {\rm Tails}  | {\rm Mon } ) =  \half. 
\end{align}
Equation \eqref{e:EPI} follows from a ``highly restricted" principle of indifference.
To justify \eqref{e:EPEI}, since SB is always woken on Monday morning the coin toss can be deferred until 
Monday afternoon. If SB is told at midday on Monday that it is Monday, 
then the (fair) coin toss is still in the future, and 
so SB should assign probability $\half$ to Heads. (We can ask SB to toss the coin.) 
Then
\be  \label{e:Elga}
 \bP(\Mon\; \&  \; \Tails) = 
 \bP( {\rm Mon } |  {\rm Tails } )  \bP(  {\rm Tails} ) =
\bP( \Tails | \Mon) \bP( \Mon),
\ee
so that  $ \bP( \Mon) = \bP(  {\rm Tails} ) =1 -  \bP( \Heads)$. 
As $\bP({\rm Heads}  | {\rm Tue} ) =0$ 
we have 
$$ \bP(  {\rm Heads} ) =  \bP( {\rm Heads}  | {\rm Mon } ) \bP(  {\rm Mon} ) 
= \half  \bP(  {\rm Mon} ) = \half - \half \bP( \Heads), $$
and so $\bP( {\rm Heads} )  = \third$. 

This argument is simple enough that it seems surprising that the conclusion should be disputed.
Nevertheless Standard Halfers reject \eqref{e:EPEI}, while Double Halfers
need to redefine what is meant by the conditioning in \eqref{e:EPEI}.

\ms 
The starting point for this paper is to ask what probability space the calculations above take place in.
In the standard SB model there is only one randomization, the toss of a fair coin,
which suggests that it should be possible to define the model on the space $\Omega_O = \{H, T\}$
with $\bP( \{H\} ) = \bP( \{T\})=\half$. 
However this cannot be the space used in Elga's calculation: no event in $\Omega_O$
has probability $\third$, and the events $\{\rm Mon\}$ and $\{\rm Tue\}$ do not lie in this space.
These events also appear to be ambiguous -- an outside observer of the experiment will experience both
Monday and Tuesday.

The set of days is  $\Omega_A = \{ \rm Mon, \rm Tue \}$, 
and so the state of SB on waking is described by the pair
$$ \om = (\om_o, \om_a ) \in   \Om_O \times \Om_A.$$
In the terminology of the philosophy literature 
elements of $\Om_O$ describe `possible worlds', while elements of $\Om_A$ give the observer's 
location in the world. The general  tendency has been to regard these two components 
as being at fundamentally different levels, a point most clearly indicated by the introduction
in  \cite{Ha05, M08, B10} of a new form of conditioning which acts in different ways 
on the two components of $\Om$.

However, when SB is woken, she is uncertain both about $\om_o$ and $\om_a$,
and it seems natural to place her ignorance about these at the same level. 
In this paper we will define our model on the space $\Om =  \Om_O \times \Om_A$.
It is viewed by two observers. The classical observer CO  only sees the space 
$\Om_O$ with its objective probability $\bP_O$. The AO sees the space
$\Om$, and the problem for the AO is to choose a suitable probability $P_A$ on $\Om$.
This is not an `extension problem' since we will not in general require that
\be  \label{e:PP0}
P_A(F \times \Om_A) = \bP_O(F) \q \hbox{  for $F \subset \Om_O$}.
\ee 
As we cannot expect to derive $P_A$ from $\bP_O$
using the axioms of probability our task divides into two parts. 

The first is to formulate reasonable properties (or `Principles') which restrict the 
possible class of measures $\PA$;
these Principles have to come from our real world intuitions.
In mathematical terms these Principles provide links between the  space $(\Om_O, \bP_O)$ 
seen by the CO and the space $(\Om, P_A)$ seen by the AO.
The second problem is to identify the class of measures allowed
by these Principles, and work out their properties. 

One Principle  mentioned in the 
philosophy literature is Lewis' Principal Principle (see \cite{L}) which roughly states that 
SB should use a probability $\PA$ which satisfies \eqref{e:PP0}.  
We will see below strong reasons for rejecting a straightforward use of this Principal Principle.
 
 The remainder of this paper is as follows. 
In Section \ref{s:MI} we introduce a generalised Sleeping Beauty problem (GSB),
where the set of occupied cells is a random subset $\sX$ of finite set $K$.
Since we are interested in inferences by the AO we introduce an auxiliary process $Z$,
and we may allow the AO in cell $x$ to observe the value $Z_x$. 
We set up the extended space $\Omega$ and formulate four Principles which give
apparently desirable properties of the AO's probability $P_A$. 
The first, denoted (PN), is a mild condition on null sets, and will be assumed throughout this paper.
The second, the  {\em Principle of Indifference} (PI), is an extension of  \eqref{e:EPI}. 
The third, the {\em Principle of Equivalent Information} (PEI), generalises \eqref{e:EPEI}. This principle
states roughly that a Classical Observer and an Anthropic Observer who are in communication and have
shared all the information they have available
should assign the same conditional probabilities to events in the space $\Omega_O$. 
The final Principle, denoted (PP), is an extension of \eqref{e:PP0}. 
We find that, except  when the number of AO is deterministic (conditional on there being any AO at all), 
there  is no probability $\PA$ which satisfies all four Principles -- see Corollary  \ref{C:4P}.
If we drop (PP) then there is a unique probability, which we denote $P_E$, which 
satisfies the remaining three Principles. 
If we require (PP) then there exists a natural probability $P_L$ which 
satisfies (PI). 

In Example \ref{E:SB} we look at the standard Sleeping Beauty problem and show that
$P_E(\Heads)=\fract13$ and $P_L(\Heads) = \fract12$; we will therefore 
refer to $P_E$ and $P_L$ as the {\em Thirder}  and {\em Halfer} measures respectively.
In Section \ref{ss:FI} (see Theorem \ref{T:PNFI}) 
we consider sequential experiments such as the original SB experiment, 
 and show that
(PEI) is equivalent to the assertion that the AO 
should assign the objective probability to random events which are still in its future. 

Section \ref{s:Obs} looks at an AO's inferences based on observation of the process $Z$. 
Section \ref{ss:ICE} shows that one needs to be  careful when one conditions on an
AO's observations, and 
Section \ref{ss:HS} looks at the model in \cite{HS}, a paper 
where insufficient  care was taken.
Section \ref{s:Ex} looks at three simple examples from astronomy and cosmology,
and we compare
the inferences made by AO which use $P_E$ and $P_L$. In the third example
we see that in a system with a large number of observers, and many independent 
components, inferences using $P_E$ and $P_L$ do not differ by much. 
(These conditions do typically hold in multiverse models.)

The very extensive  literature contains many ideas, examples and calculations.
In some cases a space such as our extension $ \Omega = \Omega_O  \times  \Omega_A$
is implicit, but the framework given here does not seem to have been set out precisely before, and 
the mathematical  formulation of the key principle (PEI) seems to be new. 
\cite{G19} does consider probability laws on 
the sample space $\Om= \Om_O \times \Om_A$, 
but rejects, at least for some probability laws on $\Om$, the possibility of conditioning
on events  like $\{\Mon\}$. Thus (PEI) would be regarded as inadmissible in the setup given there.
The calculations in \cite{R09} use an assertion similar to (PEI).

The approach of this paper, which 
is natural from a probabilistic viewpoint, concentrates on information and conditional
probabilities, and owes much to Elga's paper \cite{E}. 
An alternative approach, as seen in for example Bostrom's book \cite{B02}, is to 
assume that the observer's location has been randomly selected in some fashion
from the set of all possible locations. The type of random selection then determines 
the outcome: 
Bostrom's {\em Self Sampling Assumption} (SSA)  leads to Halfer type conclusions, while his
{\em Self Indication Assumption} (SIA) gives Thirder answers. 

Two other approaches to determining probabilities are using betting arguments or long run
frequencies. Neither of these give a clear answer for the SB problem. In general betting arguments
have to allow for the possibility of multiple bets by SB, and the question then arises how she
should regard gains or losses by her other copies. See \cite{A17} for the type of difficulties 
which arise. Given a series of (independent) SB experiments, once a week over many weeks,
it is clear that the long run proportion of `Heads weeks', that is weeks with a Heads toss,
 will be $\frac12$,  while the proportion of awakenings
in `Heads weeks' will be $\frac13$. The difficulty is how to match these limits with SB's credences
in a mathematically precise fashion. 

One reason  philosophers do not agree on the solution to the Sleeping Beauty problem
is that all known approaches lead to paradoxical conclusions.  For the Halfers 
the problems arise quickly: 
standard Halfers find themselves holding that SB's position gives her
information about the future, while Double Halfers have to redefine conditional probabilities. 
See Example \ref{E:SB} for details. 
The difficulties for Thirders do not arise  so immediately. However Bostrom's example 
of the Presumptuous Philosopher (see  \cite[p. 124]{B02} and Example \ref{E:PP}) suggests  
that the methods of inference used by Thirders, together with some  plausible extensions, 
lead to the conclusion that, without needing to make any observations at all,  
we should assign a credence of 1 to the event that the universe is infinite. 

\section{The generalized Sleeping Beauty model}  \label{s:MI}

\subsection{Probability spaces and observers} 

We generalize the Sleeping Beauty problem as follows.
Let $K$ be a finite set with $|K|=M$ where $M\ge 1$.
Let $\sP(K)$ be the set of all subsets of $K$, and let 
 $(\Omega_O, \sF_O, \bP_O)$ be a probability space carrying a random variable 
 $\sX^O : \Omega_O \to \sP(K)$,
 and a process $Z^O=(Z^O_x, x \in K)$ taking values in a finite set $\sZ$.
We assume that $0\le \bP_O( \sX = \emptyset )<1$. 
We call $(\Omega_O, \sF_O, \bP_O)$  the {\em objective probability space}. 
Most of our arguments generalize  to the case when $K$ is an infinite set 
and $\bE(|\sX^O|)<\infty$, but 
since the main challenges of the problem are already present in the
case when $K$ is finite, we will restrict to that case.

Associated with each $x \in K$ is a cell $C_x$.
The cells are identical, and are 
labelled, but on the outside, so that an occupant does not know its cell label.  
Identical anthropic observers are placed in each cell $C_x$ with $x \in \sX$.
Each cell is equipped with a telephone which permits only incoming calls. 
A Classical Observer (CO) is situated outside the collection of cells, and does not initially
know the value of $\sX^O$. The CO has a telephone which can be used to call cell 
$C_x$. By calling cell $C_x$ the CO is able to verify
whether or not $x \in \sX^O$. (We assume that the AO always answers the telephone.) 

If we wish to implement the model with a single AO, as in the original SB problem,
we choose a bijection  $\pi$ between $\{1, \dots, M\}$ and $K$. The experiment 
is run for $M$ days. On each day $j$, for $1\le j\le M$,  the AO is placed in cell
$C_{\pi(j)}$ for that day. If $\pi(j) \in \sX$ then the AO is woken on day $j$,
and at the end of day $j$ the AO's memory of that day is erased.
If $\pi(j) \not \in \sX$ then the AO is not woken on day $j$.

To handle the experience of the AO we extend the probability space as follows.
Let $\pd$ be a point not in $K$, and let
$$ \Omega_A = K \cup \{\pd \}, \q \Om = \Om_O \times \Om_A. $$
Let $S: \Omega \to \Om_A$  be defined by $S((\omega, x)) =x$ for $\om\in \Omega_O$, $x \in \Omega_A$.
Set
\be
 \FOb =\{ F \times \Omega_A, F \in \sF_O \}, \q  
 \sF =\sigma(\FOb, S) = \sigma(F \times \{x\}, F \in \sF_O, x \in  \Om_A).
 \ee
It is straightforward to check that $F \in \sF$ if and only if there exist $F_x \in \FOb$ such that
 \be  \label{e:sFrep}
   F = \bigcup_{x \in \Om_A} \{ S= x \} \cap F_x. 
 \ee

We call $\FOb$ the  {\em objective $\sigma$-field}, and events in $\FOb$ {\em objective} events.
In the terminology in the philosophy literature, events in $\FOb$ tell us about possible Worlds.
We extend $\bP_O$ to a probability measure $\bP$ on $(\Omega, \FOb)$
by setting $\bP( F \times \Omega_A) = \bP_O(F)$. 
Throughout this paper 
we will use a superscript  `$O$' to denote random variables on the base space $\Omega_O$, 
and remove it to denote the extension of the random variable  to $\Omega$.
Thus if $Y^O$ is defined on $\Om_O$, then  $Y((\om,x)) = Y^O(\om)$ for $(\om,x) \in \Om$. 

The space $(\Om, \FOb, \bP)$ describes the experience of the classical observer. Over the
course of the experiment the CO will observe that a variety of cells are occupied 
on different days, but as none of these has a  distinguished status 
the random variable $S$ is not observed by, and makes no sense to, the CO.

We now look at the viewpoint of an AO, which wakes in an initially unknown cell. 
If $G \in \sF_O$ and $x \in K$ then the interpretation of the event 
$G \times \{x\} \in \sF$ is that the  event $G$ occurs, 
and the AO is situated in cell $C_x$. 
As $S((\om,x))=x$ we see that the random variable $S$ gives the location of the AO. 
(The point $\pd$ is needed to give the value of $S$ when $\sX=\emptyset$, 
and there are no anthropic observers.)

Many papers in the probability  literature look at 
models where  observers have different information or use different probability measures. 
For example, the `expansion of filtrations' literature, following \cite{B78, J78},
looks at stochastic processes viewed by two observers with different information sets.
Systems with multiple probability measures are common in the mathematical finance literature.

The specific setup here, with the observers CO and AO, seems to be new. 
The CO uses the probability space $(\Om, \sF_W, \bP)$, while the AO observes measurable
functions on the space $(\Om,\sF)$. The question for the AO is what probability
to place on this space. 
In the theory of Markov chains on a countable state space $I$
one can derive the law $\bP^x$ of the
process started at any point $x\in I$ by conditioning from a single probability $\bP^\mu$
where  $\mu$ is a probability measure on $I$ with full support.
Similarly a law $\PA$ on $(\Om, \sF)$ 
can cover all possible locations of the AO, since the laws
$\PA( \cdot | S=x)$ give the credences of the AO if it learns that it is in cell $C_x$.

 \begin{definition} \label{D:4P}
 {\rm We define the following Principles. 
 
\sm{\bf Principle of  null sets} (PN). \\
(a)  If $F \in \FOb$ and $\bP(F)=0$ then $\PA(F)=0$. \\
(b) $P_A(S \in \sX) =1$. (It follows that $ P_A(\sX = \emptyset) = 0$.)

 \sm {\bf Principle of Indifference} (PI).  For $x \in K$, $B \subset K$, $B \neq \emptyset$, with 
$P_A( \sX = B) >0$,
 \be  \label{e:PI}
   \PA ( S=x | \; \sX=B ) = |B|^{-1} 1_B(x). 
 \ee
 
\sm {\bf Principle of Equivalent Information} (PEI). If $F \in \FOb$ and $P_A(S=x)>0$ then 
\be \label{e:PEI}
  \PA ( F  | S=x ) = \bP (F | x \in \sX ). 
\ee

\sm{\bf Principal Principle} (PP).  For $F \in \sF_W$ we have
 \be \label{e:PP1}
   \PA(F) = \bP(F| \sX \neq \emptyset).
 \ee

}\end{definition} 

The first of these places mild restrictions on null sets,  
while other three are extensions of  \eqref{e:EPI}, \eqref{e:EPEI} and \eqref{e:PP0}.

\sms
The intuition for (PI)  is that if an AO is told that $\sX= B$
and has no other information,
then it should decide it is equally likely to be in each of  cells $\{C_y, y \in B\}$. 
It is difficult to see any reason to choose an alternative distribution.

\sms
To see the motivation for (PEI), we fix $x \in K$ before the start of the experiment; this
value is known to the AO.
The CO phones cell $C_x$. If the phone is answered 
than the CO knows that the event $\{x \in \sX\}$ occurs.  If the AO receives a phone call, it
knows that it is in cell $x$. The AO and CO are in contact, and are allowed to share any 
information they have on the experiment. 
(In fact neither has anything to add to what the other already knows.)
Since they have the same information they should have identical views on the 
probability of any $\FOb$-measurable event $F$, and hence \eqref{e:PEI} should hold. 
We remark that the point of the phone conversation is just to emphasize the fact that the CO and
AO have the same information, and that it may still be reasonable to assume (PEI) in contexts
where such a conversation cannot occur. 
Note that both  \cite{Ha05, Pit15} reject assertions similar to (PEI). 

The restriction in (PEI)  that $F \in \FOb$ is necessary, since otherwise the right hand side is undefined. 
While the events $\{ S=x\}$ are disjoint, the events $\{x \in \sX\}$ are not, and 
this asymmetry may cause one to ask whether (PEI) gives the right 
mathematical description of the equivalent information of the CO and AO. However
we will see below that for some models there are strong reasons for accepting (PEI).

\begin{remark}
{\rm 
1.  None of these principles  involve the bijection $\pi$ between $K$ and 
$\{1, \dots, M\}$. 
If we had an AI which could make Bayesian inferences at a level comparable 
with the most able humans  we could consider a simultaneous experiment, where 
identical AO/AI are placed in each cell $C_x$, for $x\in \sX$.
 Some quibbles about the original SB experiment can be overcome if we use AI. 
In particular the cells can be virtual, the observations of the AO can be tightly controlled, 
and memories can be completely erased.  \\
2. 
It is conceivable that there are two different real world experiments $R_1$ and $R_2$
which are described by the same mathematical 
model (i.e. the same set $K$ and processes $\sX$ and $Z$), but where one would
have different levels of confidence in the applicability of the various Principles. 
For example, it might be very reasonable to assume (PEI) for $R_1$ but not for $R_2$.
}\end{remark}

For  $B \subset K$,  $x \in K$  define
\begin{align} \label{e:qbs-def}
  q_B = \bP( \sX = B), \q 
  Q_x = \bP( x \in \sX) = \sum_{A: x \in A} q_A.
  \end{align} 

\sms
We define a graph structure on $K$ by defining $\{x,y\}$ to be an edge if
$\bP( \{x,y\} \subset \sX)>0$. Let $E$ be the set of edges.
We set
\be \label{e:Xdef}
  X = | \sX|. 
\ee

\begin{theorem} \label{T:msets} 
(a) Suppose that $(K,E)$ is connected. 
There is a unique  probability $P_E$ on
$(\Omega , \sF)$ which satisfies (PN), (PI) and (PEI). Writing 
$\lam = 1/\bE(X)$ we have
\begin{align} 
 \label{e:FxP}
P_E( F \cap \{S=x \}) 
 &= \lam \bP( F  \cap \{x \in  \sX \}) \hbox{ for } F \in \FOb, x \in K. 
\end{align} 
In particular  for  $x \in K$,  $B \subset K$, 
\be  \label{e:st-sets}
  P_E(S=x) =   \lam Q_x, \q   P_E( \sX = B)  = \lam |B|  q_B,
\ee
and if $F \in \FOb$ then
\be  \label{e:PEF}
  P_E(F) = \frac{ \bE( 1_F X)}{\bE(X)}. 
\ee
(b) The probability $P_E$ defined by \eqref{e:FxP} satisfies (PN), (PI) and (PEI). 
\end{theorem}

\proof (a) We use the same basic strategy as in Elga \cite{E}.  
Note that the conditions on $K$ and $\sX$ imply that $Q_x>0$ for each $x \in K$.
Let $P_E$ be a probability which satisfies the hypotheses, and write
\begin{align}  \label{e:stb}
  t_B = \PA( \sX = B), \hbox{ for } B \subset K,   \q s_x = \PA( S=x )  \hbox{ for $x \in K$. }
  \end{align} 

Let $B \subset K$ with $q_B>0$ and let $x \in B$. If $t_B>0$ then by (PI)
\be \label{e:tB1}
  P_E( S=x, \sX=B) = P_E(S=x | \sX=B) t_B = \frac{t_B}{|B|} >0,
\ee
and thus $s_x>0$. If $s_x>0$ then by (PEI)
\be \label{e:sx1}
   P_E( S=x, \sX=B) = P_E(\sX=B| S=x) s_x = \bP( \sX = B | x \in \sX) s_x = \frac{ s_x q_B}{Q_x} >0,
\ee
so that $t_B>0$. 
Hence
\be
  \frac{ t_B}{q_B |B|} = \frac{ s_x }{Q_x}
\ee
for any pair $(x,B)$ with $q_B>0$ and $x \in B$. If  $\{x,y\}$ is an edge 
in the graph $(K,E)$ then it follows that $s_x/Q_x= s_y/Q_y$. As the graph $(K,E)$ is connected,
the function $s_x/Q_x$ is equal to a constant $\lam$  on $K$. We then
have  $t_B = \lam |B| q_B$ and since $\sum_{B \subset K} t_B =1$ it follows that
$\lam^{-1} = \bE( |\sX|)$.   

A further application of (PEI) proves \eqref{e:FxP}, which gives the uniqueness of $P_E$.
The final assertions \eqref{e:st-sets} and \eqref{e:PEF} follow easily from  \eqref{e:FxP}. \\
(b) Now suppose that $P_A$ satisfies \eqref{e:FxP}. 
Then $P_A(S=x) = \lam \bE( 1_\sX(x))$, and summing over $x \in K$ it follows that
$P_A( S \in K)=1$, so that $P_A( S = \pd)=0$. If $F \in \FOb$ with $\bP(F)=0$ then
$ P_A(F \cap \{S =x\})=0$ by \eqref{e:FxP}. Hence 
$$ P_A( F) = P_A( F \cap \{ S \in K\}) + P_A( F \cap \{ S = \pd\}) =0, $$
so that part (a) of (PN) holds. Finally, if $x \in K$ then by \eqref{e:FxP}
$$ P_A( S=x, x \not\in \sX ) = \lam \bP( \{ x \not\in \sX \} \cap \{ x \in \sX\})=0, $$
and summing over $x$ it follows that $P_A( S \in \sX)=1$. Thus $P_A$ satisfies (PN).

If $ B \subset K$ and $x \in B$ then by  \eqref{e:FxP}
$$ P_A( \{\sX = B \} \cap \{ S=x \} )= \lam \bP( \sX =B). $$
Summing over $x \in B$ we have 
$P_A( \sX = B ) = \lam |B|  \bP( \sX =B)$, and (PI) follows. 
Finally (PEI) follows easily from  \eqref{e:FxP}. \qed

\begin{remark} {\rm 
1. \eqref{e:st-sets} implies that 
 the $P_E$-law of $|\sX|$ is the size-biased distribution
associated with the law of $|\sX|$. This reweighting of the law of  $|\sX|$ is well known 
-- see for example \cite[p. 122]{B02}.  \\
2. If $(K,E)$ is not connected the measure $P_E$ given by \eqref{e:FxP} 
is no longer the unique measure which satisfies (PN), (PI), and (PEI). 
See the appendix to \cite{Ba}  (an earlier version of this paper) for more details on this. 
However this case does not
seem to be important in applications: all the examples that the author
has seen in the literature, when put into the framework of this paper, have $(K,E)$ 
connected.
} \end{remark}

The possibility that there are no observers does not affect the probability $P_E$.  Let 
\be \label{e:wtbP}
 \wt \bP(F ) = \bP( F |  \sX \neq  \emptyset), \q F \in \FOb,
\ee 
and let $\wt P_E$ be the unique probability which satisfies (PN), (PI), (PEI) with respect
to $\wt \bP$.

\begin{cor}
Let  $\wt P_E$ be as above. Then $\wt P_E = P_E$.
\end{cor}

\proof 
Let $p=\bP( \sX \neq  \emptyset)$; by hypothesis $p>0$. 
Let $\wt \lam = \wt\bE( |\sX|)$; then $\wt \lam = \lam/p$. By \eqref{e:FxP}
we have for $F \in \FOb$
\begin{align*}
 \wt P_E(F  \cap \{S=x\}) &= \wt \lam \wt \bP( F \cap \{x \in \sX \})  \\
&=  \lam p^{-1} \bP( F \cap \{x \in \sX \}| \sX \neq \emptyset)  \\
&= \lam  \bP( F \cap \{x \in \sX \} \cap \{\sX \neq \emptyset\} )
= \lam  \bP( F \cap \{x \in \sX \}). 
\end{align*} 
Hence $\wt P_E(F  \cap \{S=x\}) = P_E(F  \cap \{S=x\})$, and it follows that
$\wt P_E=P_E$. \qed

\sms
Let $\sX'$ be a random subset of $K$ with the property that $\sX' \subset \sX$.
Suppose that an AO learns that it is located in $\sX'$, i.e. 
that $S \in \sX'$, and wishes to estimate the probability of an event $F$.
One natural procedure is to look at
$ P_E ( F | S \in  \sX')$, but another is to look at the model associated with the 
random process $\sX'$, and then use the probability $P'_E$ obtained from that model.
More precisely, we wish to define  $P'_E$
on the space $\Om = \Om_O \times \Om_A$ with $\Om_A = K \cup \{ \pd\}$.
The graph structure $E'$ associated with $\sX'$ is given by taking $\{x,y\} \in E'$
if $\bP( \{x,y\}  \subset \sX')>0$. Even if the original graph $(K,E)$ is connected, 
the graph $(K,E')$ may not be connected, so that we cannot use Theorem \ref{T:msets}
to define $P'_E$. However, we can still define $P'_E$ by using \eqref{e:FxP}:
\be \label{e:P'Edef}
  P'_E( F \cap \{S=x \}) 
 = \lam' \bP( F  \cap \{x \in  \sX' \}) \hbox{ for } F \in \FOb, x \in K. 
\ee
Here $\lam' = 1/\bE(X')$, where $X'= |\sX'|$.

\begin{prop} \label{P:Res}
Let $\sX' \subset \sX$ with $\bP( \sX'= \emptyset)<1$, let $P_E$ be the probability defined by \eqref{e:FxP},
and $P'_E$ be given by \eqref{e:P'Edef}. Then $P_E( S \in \sX')>0$ and for $G \in \sF$, 
\be \label{e:ResP}
  P_E( G | S \in \sX') = P'_E(  G ) . 
\ee
\end{prop}

\proof It is enough to prove that for $F \in \FOb$, $x \in K$,
\be   \label{e:topr}
P_E( F \cap \{S =x \} | S \in \sX') = P'_E(  F \cap \{S=x \}).
\ee
By \eqref{e:P'Edef} the right hand side of \eqref{e:topr} is
$ \lam' \bP( F  \cap \{x \in  \sX' \}). $
For the left side we have
\begin{align*}
 P_E( F \cap \{S =x \} \cap\{S \in \sX'\})  
   &= P_E( F \cap \{ x \in \sX'\} | S=x) P_E( S=x) \\
   &= \lam \bP(  F \cap \{ x \in \sX'\}| x \in \sX )\bP( x \in \sX)  \\
   &=   \lam \bP(  F \cap \{ x \in \sX'\}).
\end{align*}
Taking $F =\Om$ and summing over $x$,
$$ P_E( S \in \sX') = \lam \sum_{x \in K} \bP(  \{ x \in \sX'\}) = \frac{ \bE(X')}{\bE(X)}. $$
Thus $P_E(S \in \sX')>0$ and 
$$  P_E( F \cap \{S =x \} | S \in \sX' )  = \lam'  \bP(  F \cap \{ x \in \sX'\}), $$
proving \eqref{e:topr}. \qed

\sms
We will call \eqref{e:ResP} the {\em restriction property} of $P_E$. 

\sms
The next result shows that it is only in rather special circumstances that we can 
find a probability which satisfies all four of the Principles given in Definition \ref{D:4P}.

\begin{cor} \label{C:4P}
(a) Suppose that there exists a law $\PA$ which satisfies (PN), (PI), (PEI) and (PP).
Then there exists $k$ such that $ \bP(|\sX| =k | \sX\neq \emptyset) =1$.  \\
(b) If there exists $k$ such that $ \bP(|\sX| =k | \sX\neq \emptyset) =1$ then
$P_L=P_E$.
\end{cor} 

\proof (a) Let $\wt \bP$  be as in \eqref{e:wtbP}. As $P_A$ satisfies
(PN), (PI) and (PEI) it is equal to $P_E$. Let $F \in \FOb$; 
using (PP) and \eqref{e:PEF} we have
\be \label{e:4pe}
 \wt \bP(F) = P_E(F) = \frac{ \wt \bE( 1_F X) }{\wt\bE(X)} \hbox{ for all } F \in \FOb.
\ee 
Let $k$ be such that $\bP(X=k)>0$. Setting $F =\{X=k\}$ and using \eqref{e:4pe} 
gives that $\wt \bE (X)=k$, so that there is exactly one such $k$.  \\
(b) is immediate from \eqref{e:PEF}. \qed

If one wishes to keep (PP) 
then the most natural probability to use is given by
 \be \label{e:PL}
  P_L( F  \cap \{ \sX = B \} \cap \{ S=x\} ) = 
  \begin{cases}
  |B|^{-1} \wt \bP(F \cap \{\sX=B\} )  \hbox{ for } x \in B, B \subset K, \\
  0 \hbox{ for } x \notin B,  B \subset K.
  \end{cases}
 \ee
This satisfies  (PN), (PI) and (PP).  We can 
strengthen (PI) by requiring that, for $x \in  B$, $B \subset K$, and $F \in \sF_W$, 
\be   \label{e:PIst} 
   \PA ( S=x | \; F \cap \{\sX=B \} ) = |B|^{-1} \q \hbox{ if }   \PA(F \cap \{\sX=B\})>0.  
 \ee

\begin{prop} \label{P:PL} 
The probability $P_L$ is the  unique  probability on 
$(\Omega , \sF)$ which satisfies (PN), (PP) and \eqref{e:PIst}.
\end{prop}

\proof Let $\PA$ satisfy (PN), (PP) and \eqref{e:PIst}. Then if $x$, $B$ and $F$ are as above
and $x \in B$, 
\be  \label{e:pl1}
   \PA ( \{ S=x \} \cap F \cap \{\sX=B \} ) = |B|^{-1} \bP(  F \cap \{\sX=B \} | X \ge 1). 
\ee   
If $x \not\in B$ then (PN) implies that the left side of \eqref{e:pl1} is zero. Thus $\PA$
satisfies \eqref{e:PL} and so equals $P_L$. \qed

Easy examples show that $P_L$ does not satisfy the  restriction property.
(See also Section \ref{ss:2phase}.)

\begin{example} \label{E:SB}
{\rm {\em The original Sleeping Beauty problem.}   
Let $K=\{1,2\}$, and $\Om_O =\{0,1\}$, with $\bP_O(\{0\})=\bP_O(\{1\})=\half$, and
let $\sX^O(0)=\{1\}$, $\sX^O(1)=K$. Since $\bP_O(\sX^O=\emptyset)=0$ we can dispense with
the extra point $\pd$, and take $\Om_A = K$, so that $\Om = \Om_O \times K$.
As usual we write $\bP$ and $\sX$ for the extensions of $\bP_O$ and $\sX^O$ to $\Om$.

Let $\PA$ be a probability on $\Om$ which satisfies  (PN), so that
$\PA( \{ (0,2) \})=0$.
Set
\be \label{e:SB3prob}
 a= \PA(\{(0,1)\}), \, b= \PA(\{(1,1)\}), \, 1-a-b = \PA(\{(1,2)\}).
\ee
We define the events $\Heads=\{ (0,1), (0,2)\}$, $\Mon=\{ (0,1), (1,1)\}$, and let
$\Tails$ and $\Tue$ be the complements of $\Heads$ and $\Mon$ respectively.
Then
\be \label{e:condP1}
 \PA(\Heads | \Mon) = \frac{a}{ a+b}, \hbox{ and  }
\PA(\Mon|\Tails) = \frac{b}{1-a}.
\ee
(PEI) states that $\PA(\Heads | \Mon) = \half$, which implies that $a=b$.
(PI)  gives $\PA(\Mon|\Tails) =\half$, so that $a+2b=1$. Finally
(PP) gives $a=\half$. So each of these principles imposes a different 
linear relation between $a$ and $b$, and thus the three conditions are independent. 

Taking (PI) and (PEI) together we obtain $a=b=\frac13$, and so
we have $P_E(\Heads)=\frac13$. 
If (PI) and (PP) hold then $a=\half$, $b = \fract14$, and  $\PA=P_L$. 
We have  $P_L(\Heads | \Mon) = \frac23$, so this measure does not satisfy (PEI).  

The interpretation of this conditional probability is a difficulty for Halfers. 
Even though the coin toss can take place
on Monday afternoon, on the standard Halfer view SB at midday on Monday should hold 
that this fair coin toss has probability $\fract23$ of coming down Heads.
Lewis \cite{L-E} did adopt this position, but nearly all subsequent authors have found
it unsatisfactory. 
Double Halfers hold that both the probability of Heads, 
and  the probability of Heads given that it is Monday are $\half$. 
If we use  \eqref{e:SB3prob} and \eqref{e:condP1} then we obtain 
$a=\half$ and $a/(a+b)= \half$, so that $a=b=\frac12$ and $P_A(\Tue)=0$. 
Thus Double Halfers would appear to be committed to the view that 
on waking SB is certain that it is Monday: \cite{Ha13} does adopt that position. 
But most Double Halfers wish to allow $\PA(\Tue)>0$, and some
wish to keep (PI) and hold that $\PA(\Mon|\Tails) = \half$. (Thus they are really `Triple Halfers'). 
It is clear from the calculations above that
this position cannot be described by the usual axioms of probability. 
Double Halfers are well aware of this difficulty, and the literature
contains a number of strategies for dealing with it. As one might expect, attempts 
such as those in  \cite{Ha05, M08, B10}
to modify the standard definition of conditional probability leads to operations with 
mathematically undesirable properties.  See  \cite[Sections 4, A.2]{Ba} for more details.
}\end{example}

\begin{example} \label{Ex:Pittard}
{\rm {\em Four Beauties.} 
This elegant example is due to Pittard \cite{Pit15}.
The experiment lasts for one night only, with 
four participants (`Beauties') A, B, C, D. One of the four will be chosen at random to be the  `Waker'.
The four participants all go to sleep. Three times in the night a pair of participants are woken
at the same time for a period (say 15 minutes) during which they can converse. They then
go back to sleep, and as usual are given a memory erasing drug so that they forget the awakening.
The Waker is woken all three times, and paired with each other participant exactly once.

We analyse this model by describing the awakenings of a particular 
participant, A say. We take  $\Om_O =  \{A, B,C,D\}$, and let
$K= \Om_A= \{ B,C,D\}$ be the set of possible partners of A during an awakening. 
We write $\om=(v, x)$ for points in $\Om$, 
and with $\om$ of this form we define $W(\om) = v$, $S(\om) = x \in K$.
So $W$ gives the Waker, and $S$ describes A's partner for the selected awakening.
$\bP$ is the probability which makes $W$ uniform on $\Om_O$, and
$\sX$ is defined by $\sX = \{K\}$ if $W=A$, and $\sX= \{W\}$ if $W \neq A$.

As $\bP(\sX = K)>0$ the model satisfies the hypotheses  of Theorem  \ref{T:msets}.
We have $Q_x = \half$ for $x=B,C,D$, and
$\bE| \sX| = \frac32$, so that $\lam=\frac23$.
Thus $P_E(S = x )= \lam Q_x = \frac13$ for $x=B,C,D$.
Suppose that A wakes and sees that the other person woken is B.
Then, using (PEI)  the probability that A is the Waker is
$$ P_E( W=A | S=B ) = \bP( W=A |  B \in \sX)
= \frac{ \bP( W=A)}{ \bP(  B \in \sX) }=\fract12.   $$
This answer accords with common sense --
A and B are in the same situation and by symmetry each regards herself as equally likely to be the Waker.

On the other hand, if we use $P_L$ then $P_L(W=x) = \bP(W=x) = \fract14$,
for each $x \in K$. By (PI) $P_L(S=B| W=A) = \fract13$, and so
$P_L(S=B, W=A) = \fract{1}{12}$. 
By symmetry $P_L(S=B) = \fract13$, 
%
and hence
$$ P_L( W=A| S=B) = \frac{ P_L( S =B , W=A) }{ P_L(S=B) } =  \fract14. $$
Thus if A and B are woken together then as A has credence $\frac14$ that she is the Waker,
she must have credence $\frac34$ that B is the Waker. 
 B, of course, has the opposite credences. 
A and B appear to have the same information, and can share all they that know, but still
have different credences for an objective event.
Although this example would seem to pose a severe challenge for Halfers, \cite{Pit15} 
does nevertheless adhere to the Halfer position.
}\end{example}

\begin{example} \label{E:PP}
{\rm 
{\em The  Presumptuous Philosopher.}
(For an explanation of the terminology see \cite{B02}.) 
Here we take $1\le N \ll M$, and set $K = \{ 1, \dots, M\}$. Let $W$ be a $\Ber(\half)$ r.v
and let $\sX = K$ if $W=1$ and $\sX = \{ 1, \dots , N \}$ if $W=0$.
Using \eqref{e:PEF} and \eqref{e:PP1} we have
\be
   P_E(W=0) = \frac{N}{M+N}, \q P_L(W=0) = \half. 
\ee
If we take $|\sX|$ to be the size of the universe, then an AO which uses $P_E$ strongly
favours the larger universe. If say $N=10^{25}$ and $M=10^N$ then the $P_E$ probability
that $W=0$ is so small that any observation to the contrary is overwhelmingly likely
to be due to experimental error.

Bostrom \cite{B02} rejects the Thirder solution on the basis of this example. But,
as we have seen, the Halfer position also runs into difficulties, which seem to me to be 
more immediate and fundamental than the problem here.
}
\end{example} 
\subsection{Future information} \label{ss:FI}

In the introduction we noted that for the standard SB problem, one argument for \eqref{e:EPEI}
is that as the coin can be tossed on Monday evening, on Monday morning it is a future event, and so 
the AO and CO should assign it the same probability. 

Our model does not have a natural time structure; to implement it for a SB type 
experiment we needed to choose a bijection between the set of states
$K$ and $\{1, \dots, M\}$. However, an important subclass of models does have
such a structure; these include the original SB problem and the Doomsday
problem (see \cite{D}). For these models we have 
\be  \label{e:DAm}
 \sX = \{ k, 1\le k \le X\},
\ee
and the interpretation is that on days $1, \dots, X$ the AO is woken in a cell and then
(as usual) has  its memories erased at the end of each day. There are no awakenings
after day $X$, at least until the experiment is concluded. 
We define $q_n, Q_n$ for $0\le n \le M$ by \eqref{e:qbs-def}, and 
assume that $M$ is chosen so that $q_M = \bP( X=M)>0$.

We now show that, subject to a mild restriction on the probability space, 
for models which satisfy \eqref{e:DAm} 
(PEI) is equivalent to a condition which states that the AO and CO have the same
views on the probability of future events. In essence this argument proceeds by pushing
every randomisation not needed to determine whether or not $X\ge n$  into the post-$n$ future.

To avoid measure-theoretic technicalities, we make the assumption that the objective probability 
space $(\Omega_O, \sF_O, \bP_O)$ has the following structure. 

\begin{assumption} \label{A:omreg}
{\rm
(1) There exist independent Bernoulli random variables $(V^O_n, 0\le n \le M)$ and $U^O$, such that
$U^O$ is uniform on $[0,1]$ and 
\be
  \bP_O( V^O_n = 1 ) = 1 -  \bP_O(  V^O_n   = 0)  = q_n/Q_n, \q 0 \le n \le M.
 \ee 
(2) $\sF_O = \sigma ( U^O,  V^O_n  , 0 \le n \le M)$.\\
(3) The r.v $X^O$ satisfies
\be \label{e:xodef}
  X^O = \min\{ n:  V^O_n  =1\}. 
\ee
}\end{assumption} 

\begin{remark}  \label{R:probspace}
{\rm 
This assumption is much less restrictive than it might appear. 
A {\em standard probability space}  is a space which is isomorphic to the unit interval 
with Lebesgue measure. In  \cite[p.~61]{I84} It\^o remarks that ``all probability spaces appearing 
in practical applications are standard'', which can be restated as saying
that any stochastic process appearing in practical applications can be defined as a
measurable function of a single r.v.~$U$ with uniform distribution on $[0,1]$. 
In particular, given a probability space $(\Om', \sF', \bP')$ carrying a random variable 
$X'$ taking values in $\{0, \dots, M\}$ and a stochastic process
$Z'=(\wt Z_1, \dots , \wt Z_M)$  taking values in a complete separable metric space $\sK$,
then there exists a probability space $(\Om_O, \sF_O, \bP_O)$ satisfying 
Assumption \ref{A:omreg} carrying  $X^O$, $Z^O$ such that
$(X^O, Z^O_1 \dots , Z^O_M)$ has the same distribution as $(X', Z'_1, \dots , Z'_M)$.
See \cite{I84} for more details. 
} \end{remark}

We write $(V_n)$, $U$ for the extensions of $( V^O_n , n \in \bZ_+)$ and $U^O$ to the space
$(\Omega, \sF, \bP)$.  Assumption  \ref{A:omreg} implies that $\FOb= \sigma(U, V_n, 0\le n \le M)$.
From \eqref{e:xodef} we have
$$ \{X \ge n\} =\{ V_k=0, 0\le k \le n-1\}. $$ 
We regard the randomization $V_k$ as being made at the end of day $k$, and the randomization 
$U$ being made at the end of day $M$. 

 \begin{definition}
 {\rm 

\sm{\bf Principle of  no future information} (PNFI). 
Let $P_A$ be a probability on $(\Om, \sF)$ which satisfies (PN).
$P_A$ satisfies (PNFI) if whenever $0 \le n \le M$, and 
$G \in \sigma( V_n, \dots , V_M, U)$, then
\be \label{e:pnfi}
  \PA(G | S=n) = \bP(G).   
\ee
}\end{definition} 

We assume that the CO is able to observe the random variables $V_k$ as they occur.
The value $n$ is fixed at the start of the experiment. On 
day $n$ the CO knows the values of $V_0, \dots, V_{n-1}$,
and thus whether or not $X \ge n$. If $X<n$ the CO takes no action. 
If $V_0=\dots = V_{n-1}=0$, so that $X \ge n$, then as 
$G$ is independent of $V_0, \dots, V_{n-1}$ we have $\bP(G | X \ge n) = \bP(G)$,
and the right hand side of \eqref{e:pnfi} gives the CO's conditional probability 
that $G$ will occur.  
If $X \ge n$ then the CO telephones cell $n$ and talks to the AO there. The AO 
thus knows that the event $\{S=n\}$ holds, and hence also knows that $\{X \ge n\}$.
The event $G$ lies in the future for the AO also, and the left side of \eqref{e:pnfi} is the AO's 
conditional probability that it will occur.  (PNFI) then states that the two observers agree
on the probability of this future event.

\begin{thm} \label{T:PNFI}
Suppose that Assumption \ref{A:omreg} holds. Then (PEI) and (PNFI) are equivalent
for probabilities $P_A$ which satisfy (PN). 
\end{thm} 

\proof Suppose that (PEI) holds, and let $G \in \sigma( V_n, \dots , V_M, U)$.
Then $G \in \FOb$ and is independent of $\{X \ge n\}$.
Hence \eqref{e:pnfi} follows immediately from (PEI). 

Now suppose that (PNFI) holds. Let $(e_0, \dots , e_{n-1} ) \in \{0,1\}^n$, 
$H = \{ V_k = e_k, 0\le k \le n-1\}$, and $H_1 = \{ V_k =0,   0\le k \le n-1\}$.
Let  $G \in \sigma( V_n, \dots , V_M, U)$. If the AO knows that $\{S=n\}$ then it knows
that $H_1$ has occurred, and so $\PA( G \cap H | S=n) =0$ if $H \neq H_1$, while
\be
   \PA(  G \cap H_1 | S=n) = \PA(G| S=n).
\ee
Similarly we have $\bP( G \cap H | X \ge n) =0$ if $H \neq H_1$, and
 $\bP( G \cap H_1 | X\ge n ) =\bP(G)$. We therefore deduce that 
\be
   \PA(  G \cap H | S=n) = \bP( G \cap H |X \ge n ).
\ee
Using the $\pi$-$\lam$ Lemma (see \cite[p. 447]{Du}) and property (2)
of Assumption \ref{A:omreg} it follows that (PEI) holds. \qed

\section{Observations by Anthropic Observers }  \label{s:Obs}

\subsection{Technicolour Beauty} \label{ss:aux}

\sms
Recall that $Z=(Z_x, x \in K)$ is an $\FOb$ measurable process
taking values in a finite set $\sZ$.
We call $Z_x$ the {\em colour} of cell $C_x$ and 
assume that the AO is able to observe this value. 
This model is often called {\em Technicolour Beauty} --  see \cite{T08, CKS3}.
For $y \in \sZ$ let
\be  \label{e:Lydef}
  L_y = |\{ x \in \sX: Z_x = y \} | = \sum_{x \in K} 1_{(x  \in \sX )} 1_{ ( Z_x =y ) }  
\ee
be the number of AO which see the colour $y$. Set
\be \label{e:GHydef} 
 G_y =  \{L_y \ge 1 \}; \q  H_y = \{S \in \sX, Z_S =y\}.  
 \ee  %
Thus $G_y$ is the event that some AO sees $y$, while $H_y$
is the event that  the AO selected by $S$ sees the value $y$. 
We have $G_y \in \FOb$, and in general $H_y \not\in \FOb$. 
We also have $H_y \subset G_y$. 

\sms
Conditioning on the observation $H_y$ does not affect the size biasing 
obtained by using $P_E$. 

\begin{lemma} \label{L:ZS}
For $F \in \FOb$ and $y \in \sZ$, 
\be  \label{e:PEZ}
  P_E( F \cap \{ Z_S=y\}) = \lam \bE(1_F L_y).
\ee
Hence 
\be \label{e:ZSy}
  P_E( F | Z_S=y) = \frac{ \bE(1_F L_y) }{\bE(L_y) }. 
\ee
\end{lemma}

\proof 
We have using \eqref{e:st-sets}  and (PEI) 
\begin{align*}
   P_E(  F \cap H_y ) 
 &= \sum_{x \in K} P_E( F \cap \{ Z_x =y \}| S=x  ) P_E(S=x) \\  
   & =  \sum_{x \in K}  \bP( F \cap \{ Z_x =y \} | x \in \sX ) \lam Q_x   \\
   &=  \lam  \sum_{x \in K}  \bP( F \cap \{ Z_x =y \} \cap \{x \in \sX \} ) = \lam \bE( L_y 1_F) .
\end{align*}
Setting $F = \Omega$ gives $P_E(H_y) = \lam \bE(L_y)$, and \eqref{e:ZSy} follows. \qed

We compare this with what happens when we use $P_L$.

\begin{lemma} \label{L:ZPL}
For $F \in \FOb$ and $y \in \sZ$, 
\be \label{e:PLZ}
 P_L( F \cap \{Z_S=y\}) 
  =  \bE \Big(  \frac{ 1_F L_y}{X} \Big | X\ge 1 \Big), 
\ee
and 
\be \label{e:ZSLy}
  P_L( F | Z_S=y) 
  = \frac{ \bE\big(  { 1_F L_y}{X^{-1}} \Big | X\ge 1 \big) }
      {  \bE\big(  { L_y}{X^{-1}} \Big | X\ge 1  \big) }.
\ee
\end{lemma}

\proof We have
\begin{align*}
 P_L( F \cap \{ Z_S=y \} ) 
 &= \sum_{B \neq \emptyset} \sum_{x \in B} P_L( F \cap \{\sX=B\} \cap \{ Z_x=y\} \cap \{S=x\} ) \\
 &= \sum_{B \neq \emptyset} \sum_{x \in B} |B|^{-1} \bP( F \cap \{\sX=B\} \cap \{ Z_x=y\} | \sX \neq \emptyset) \\
  &= \sum_{B \neq \emptyset}  |B|^{-1} \bE( 1_F  1_{\{\sX=B\}} L_y \} | \sX \neq \emptyset) 
  = \bE\big(  { 1_F L_y}{X^{-1}} \Big | X >0  \big).
\end{align*}
Dividing this expression with the one for $P_L(Z_S=y)$ gives \eqref{e:ZSLy}. \qed

\begin{remark}
{\rm 
The presence of the term $X^{-1}$ in \eqref{e:PLZ} 
means one seldom obtains simple expressions for $P_L$ probabilities.
}\end{remark}


\begin{defn}
{\rm  We say that an auxiliary process $Z$ is {\em injective} if $\bP(L_y \le 1)=1$ for each $y \in \sZ$. 
 }\end{defn}

\begin{lemma}  \label{L:TBF}
Suppose that $Z$ is injective. Then 
\be  \label{e:PXG}
  P_E( F | H_y) = \bP(F| G_y) \q \hbox{ for all $F \in \FOb$.}
\ee
\end{lemma} 

\proof 
As $Z$ is injective $L_y = 1_{G_y}$ 
and \eqref{e:PXG} follows from \eqref{e:ZSy}. \qed

\subsection{Improper conditioning events and symmetry breaking}   \label{ss:ICE}

Let $(\Om, \sF, \bP)$ be a probability space, and 
$F,G,H$ be events with $H \subset G$ and  $\bP(G \setminus H )>0$. 
The conditional probability of $F$ given $H$ is 
\be
 \bP(F | H ) = \frac{ \bP(F \cap H)}{\bP(H)},
\ee
which is often informally described as being ``the  probability that $F$ will occur given 
that the observer knows that $H$ occurs". 
However if $H$ occurs then $G$ also occurs, 
so that this informal description does not explain why one could not use $\bP(F \cap G)/\bP(G)$
instead. 
This procedure is so obviously incorrect that few if any introductions to 
conditional probability consider this possibility. 
(Some works are quite careful in their explanations of conditional probability and avoid
the informal description above -- see for example  \cite[p.~115]{F57}.) 
In terms of Kolmogorov's definition of conditional expectation, one is interested in
\be \label{e:kolm}
 \bE( 1_F | \sigma(1_H)) = 1_H \bP(F|H) + 1_{H^c}  \bP(F|H^c). 
\ee
If one replaces $H$ by $G$ in \eqref{e:kolm} then the resulting random variable is not measurable
with respect to the agent's observation, which is of $1_H$.
 I will call conditioning events of the kind that $G$ is here {\em improper conditioning events}.
 Easy examples show that replacing $H$ by an improper conditioning event $G$
 can lead to a wildly incorrect answer.

\ms 
In the standard SB experiment she has identical awakenings on each day, 
but one can break this symmetry by introducing  an auxiliary process $Z$ as in Section \ref{ss:aux}
and allowing SB to observe the value of $Z$ in her cell. 
This suggests that one might be able to describe SB's experience  by using the
objective probability $\bP$. 
One paper which adopts this approach is \cite{CKS3}. Using the notation introduced
above 
\cite{CKS3} gives the definition 
\be \label{e:PZb}
   P_{Z,y}(F ) = \bP(F | G_y) \hbox { for  } F \in \FOb.  
\ee
The interpretation in \cite{CKS3} is that if the AO observes that the process 
$Z$ in its cell takes the value $y$  (i.e.~the event $H_y$) then it knows that $G_y$ occurs, 
and so \eqref{e:PZb} gives the AO's credence that $F$ occurs.
While of course the definition \eqref{e:PZb} makes mathematical sense, 
since $G_y$ is an improper conditioning event it is questionable whether the 
law $P_{Z,y}$ correctly models the experience of the AO. 
If $Z$ is injective then $P_{Z,y} = P_E(\cdot | H_y)$ by Lemma \ref{L:TBF},
but as \cite{CKS3} allows for more general $Z$ it is not surprising that 
by varying the process $Z$ they are able obtain a range of values for $P_{Z,y}(\Heads)$.
The framework of \cite{CKS3} does not allow one to combine the laws $P_{Z,y}$ into a
single probability $P_Z$ which gives SB's credences before she makes an observation.
As the authors of  \cite{CKS3} remark, they cannot use the law of total probability since the
events $G_y$ are not disjoint. 

\sms
The next section looks at an example of the use of 
 improper conditioning events in the physics literature.

\subsection{A model of Hartle and Srednicki}  \label{ss:HS}

The paper \cite{HS} 
questions arguments, such as those in \cite{PageD}, which assume that we are `typical' observers. 
To help clarify the issues they introduce a model of a universe which has $N$
successive cycles, and assume that observers in cycle $i$ have no knowledge
of the state of the universe in any other cycle. Each cycle  has one of two global
properties: red (R) or blue (B). In each cycle the probability that an `observing system' (such as
`us') exists is $p\in (0,1)$, and is independent of whether the global state is red or blue. 

There are two theories of the universe, which are equally likely. 
The first, $AR$ is that all the cycles are red, and the second,
$SR$, is that exactly $M$ of the cycles are red, and $N-M$ are blue. Here $1\le M \le N-1$. 
Hartle and Srednicki do not specify which cycles are red or blue in the $SR$ case; for simplicity
we will assume that under $SR$ each possible arrangement of $M$ red and $N-M$ blue is equally likely. 
For each theory we write $N_R(T)$ and $N_B(T)$ for the number of red and blue cycles according to the theory $T$.
(Thus $N_R(AR)=N$ and $N_R(SR) = M$.) 

The authors of \cite{HS} wish to calculate the probability of SR given that `we' (situated in one of the cycles)
observe that our cycle is red. Writing $H_{ER}$ for the event ``we exist and observe red",
they claim that
\be \label{e:hsw1}
   P(H_{ER}|T)  = 1 - (1-p)^{N_R(T)},
\ee
and then use Bayes' formula to compute $P(T|H_{ER})$. Writing $f(p,n) = 1  - (1-p)^n$, they obtain
\be  \label{e:hs-bay}
  P(SR| H_{ER}) = \frac{ f(p,M) }{ f(p,N) + f(p,N) } .
\ee
This conclusion is counterintuitive. If for example $M=1$, $N=1000$ and $p$ is very close to 1,
then \eqref{e:hs-bay} gives that $ P(SR| H_{ER}) \approx \half$: our observation of red has 
given very little information even though, on the hypothesis SR, red universes are very rare.

However \eqref{e:hsw1}, and hence  \eqref{e:hs-bay}, are incorrect. 
To see that  
\eqref{e:hsw1} cannot be true, let $H_E$ be the event ``we exist", and
note that by red/blue symmetry we should also have
$P(H_{EB}|T) = 1 - (1-p)^{N_B(T)}$. Since  $P(H_{ER}|T)+P(H_{EB}|T) = P(H_E|T) \le 1$,
we obtain
\be  \label{e:hsw2}
   2 - (1-p)^{N_R(T)} - (1-p)^{N_B(T)}  \le 1,
\ee
and  this fails for many values of $p$, $N_R(T)$ and $N_B(T)$.

The error which led to \eqref{e:hsw1} is confusion between the events 
$H_{ER} =\{$we exist and observe red$\}$ and $G_{ER} =\{$there is a cycle in which an observer exists
and observes red$\}$. 
Since $H_{ER} \subset G_{ER}$ if we observe that $H_{ER}$ occurs then we 
also observe that $G_{ER}$ occurs. The events $G_{ER}$ and $G_{EB}$ are not disjoint, which is 
why the sum on the left hand side of \eqref{e:hsw2} can be greater than 1.
I suspect that the cause of  this error is that the authors of \cite{HS} were (rightly) uneasy 
about the event $H_{ER}$, and wished to find a suitable objective event which they could work with.
In the terminology of the previous section, $G_{ER}$ is an improper conditioning event. 

We now apply our formalism to this problem. We define
$$ \Omega_O = \{AR, SR\} \times  ( \{ R, B\} \times \{0,1\} )^N . $$
An element $ \om \in \Omega_O$ is a sequence
\be \label{e:omdef} 
 \om = (y, a_1, n_1, a_2, n_2, \dots , a_N, n_N ), 
 \ee
where $y \in \{ AR,SR\}$, and for $1\le i\le N$ we have $a_i \in \{R,B\}$ and $n_i \in \{0,1\}$.
Thus $a_j$ gives the state (red or blue)  of cycle $j$, and $n_j=1$ if there are observing systems in cycle $j$,
and $n_j =0$ otherwise. 
We define $\Omega_A = \{ \pd, 1, \dots, N\}$.
With $\om$ given by \eqref{e:omdef} we define random variables on $\Omega = \Omega_O \times \Omega_A$
by
\be
  U(\om,k ) = y, \q S(\om, k) =k, \q Z_i(\om,k) = a_i, \q   \xi_i (\om,k) = n_i, \hbox{ for } 1\le i \le N  .
\ee
The random set of occupied universes is $\sX = \{ i : \xi_i =1 \}$.
Under the objective probability law $\bP$  the random variables $U$, $(\xi_j, 1\le j \le N)$ are
 independent with
$\bP(U=AR)=\bP(U=SR)=\half$, and $\bP(\xi_j =1) = 1-\bP (\xi_j =0) = p$. If $U=AR$ then 
$Z_i = R$ for all $i$, while if  $U=SR$ then  exactly $M$ of the $Z_i$ are equal to $R$, and the remainder 
are equal to $B$, with each ordering of reds and blues being equally likely. 
We assume that (PN), (PI) and (PEI) hold and
write $P_E$ for the probability which satisfies these conditions.
Note that $H_{ER} = \{Z_S = R\}$. 
For $y \in \{ R,B\}$ let
$$ L_y = | \{ j: \xi_j =1, Z_j = y \} | = | \sX \cap \{ i: Z_i =y\}|. $$
By Lemma \ref{L:ZS}, for $F\in \FOb$,
$$ P_E (F | Z_S =y) = \frac{ \bE( 1_F L_y) }{\bE(L_y) }. $$
Taking $F= \{ U =SR\}$ we have 
$\bE( 1_F L_y) = \half \bE( L_y | U = SR) = \half p M$,
and similarly $ \bE( L_y) = \half p (M+N)$. Hence
\be \label{e:MNM}
    P_E ( U= SR | Z_S = R) = \frac{ M }{ M + N}. 
\ee
In the case noted above, when $M=1$, $N=1000$, we obtain the natural and unsurprising
result that $P_E ( U= SR | Z_S = R) = 1/1001$.

By Corollary \ref{C:4P} the laws $P_E$ and $P_L$ for this model are not the same, but
using Lemma \ref{L:ZPL} we find that one does have $P_L ( U= SR | Z_S = R) = M/(M+N)$. 


\section{Examples -- Life, the Universe, and Everything } 
 \label{s:Ex}

\subsection{Probability of life } \label{ss:Plife} 

Look at $M$ sites (e.g. planets, each associated with a separate star) 
so that $K=\{1, \dots , M\}$. 
Let $V$ be a r.v. on $(0,1)$ with a continuous density function $f(v)$, and conditional on $V$ let 
$(\xi_x, x \in K)$ be independent $\Ber(V)$ r.v. 
We set $\lam = 1/\bE(V)$. We define $\sX = \{ x: \xi_x=1\}$, and $X=|\sX|$.
We interpret $V$ as the probability life develops on a planet, 
and say that the planet  $x$ has life if $\xi_x=1$. For simplicity we assume that life is the same
as the existence of observers, and neglect any consideration of how long this might take.
(See \cite{VH, SDT} for more realistic models.)

We have $\bE(X | V ) = MV$, and 
the probability that no observers exist is 
$$ p_M= \bP(X=0) = \int_0^1 (1-v)^M f(v) dv =\bE (1-V)^M . $$


We now look at the density of $V$ given a variety of possible observations. 

\sm {\em (0) An external observer with no observation.}
In this case the observer just uses the prior density $f(v)$.

\sm {\em (1) An external observer.} Suppose that this observer looks at site $x$ and observes $\xi_x=1$.
Standard Bayesian computations give
$$ \bP( V \le t| \xi_x=1) = \frac{ \bP(V \le t, \xi_x=1) }{\bP( \xi_x=1) }    
= \frac{ \int_0^t v f(v) dv }{  \int_0^1 v f(v) dv } = \lam  \int_0^t v f(v) dv. $$
So the posterior density of $V$ conditional on $\xi_x =1$ is
\be \label{e:f1}
 f_1(v) = \lam v f(v). 
\ee 
 
\sm {\em (2) An AO using $P_E$.} 
This AO is on planet $S$, and observes $\xi_S=1$.
By Theorem \ref{T:msets}
$$ P_E( V \le t )   = \frac{ \bE( X 1_{(V \le t)} )}{\bE(X) } 
= \frac{ \bE( V 1_{(V \le t)} )}{\bE(V) } = \lam  \int_0^t v f(v) dv. $$
So the AO's `posterior' density  is the same as  in (1). 

\sm {\em (3) An AO  using $P_L$.} 
By \eqref{e:PP1} 
$$ P_L( V \le t) = \bP(  V \le t | \sX\neq \emptyset) = (1-p_M)^{-1} 
\int_0^t \Big( 1 - (1-v)^M \Big) f(v) dv .  $$
In the $M \to \infty$ limit this gives the density
$f_{3,\infty}(v) = f(v)$. So in the large $M$ limit for an AO using $P_L$ the observation $\xi_S=1$
gives no information. 

\begin{remark} \label{R:M1}
{\rm 1. The case $M=1$ can be interpreted as a single universe, where the constants of physics
are  chosen according to some random distribution and giving rise to observers with probability $V$.
By Corollary \ref{C:4P} we have $P_L=P_E$, and the AO which uses either of these
will adopt the posterior density $f_1$ given by \eqref{e:f1}. So we find, somewhat remarkably,
that an AO gains information  from its existence, even though it is not possible
for it to observe that it does not exist.  \\
2. In the standard SB problem the law $P_L$ looks initially plausible (see \cite{R09}),
but is found to give unreasonable results when one conditions on a suitable event $F$.
We can loosely call these {\em challenge events} for $P_L$. 
In  this example there do not seem to be any natural challenge events for $P_L$,
and so, looking at this example on its own, it does not seem
easy to decide which of $P_E$ or $P_L$ one should adopt.  
}\end{remark}

\subsection{A two zone universe} \label{ss:2phase}

We modify the example above by considering a universe with two zones. 
Let $K= \{ -M, \dots, -1\} \cup \{ 1, \dots, M\}$, let $p_0, p_1, p_2$ be strictly positive, 
$W\sim  \Ber(\half)$, let $V_i = p_0$ for $x<0$ and $V_x = p_1 1_{(W=0)} + p_2  1_{(W=1)}$
for $x>0$.  
Conditional on $W$, let $(\eta_x, x \in K)$  be independent $\Ber(V_x)$ r.v.
We set $\sX = \{ x: \eta_x =1\}$.
We assume that $M$ is large enough so that $\bP(X=0)$ can be neglected.
Define $K_0 = K \cap \bZ_-$, $K_1 = K \cap \bZ_+$, $\sX_j = \sX \cap K_j$,
and $X_j = |\sX_j|$, $X = X_0 + X_1$. 

We have
\be
P_L(V_1=p_1) =   P_L( W = 0) = \bP( W=0| X \neq 0) \simeq \bP(W=0) = \half. 
 \ee

Now let $Z_i = 1_{(i \ge 1)}$. By Lemma \ref{L:ZPL}, and noting that
$L_0(Z) = |   \{ i \in \sX: Z_i =0 \} | = X_0$, 
$$ P_L( V_1= p_1 | Z_S =0) = \frac{ \bE( X_0 X^{-1} 1_{(W=0)} | X>0)}{ \bE( X_0 X^{-1} | X>0) }. $$
For large $M$ we have 
$X_0/M \simeq  p_0$ and $X_1/M \simeq  V_1/M$. Hence 
$$  \bE( L_0 X^{-1} 1_{(W=0)} | X>0) \simeq \half  p_0/( p_0 + p_1), $$
and
$$ P_L( V_1= p_1 | Z_S =0) \simeq 
\frac{  p_0/( p_0 + p_1) }{  p_0/( p_0 + p_1)  +  p_0/( p_0 + p_2) }  
= \frac{p_0 +p_2}{2p_0 + p_1 + p_2}. $$ 

If $p_0, p_1 \ll p_2$ then this probability is close to 1.
Thus using  the terminology
of the previous example, $\{Z_S =0\} =\{ S \in K_0\}$ is  a challenge event for $P_L$. 
Initially an observer (in some planet $S$) has the view that $P_L( V_1 = p_1) \simeq \half$. If 
further observation reveals that $S$ is in the negative zone $K_0$,
then the observer does not seem to have gained any information on 
conditions in the positive zone $K_1$, 
but  nevertheless now assigns a probability close to 1 to the event $\{V_1 = p_1\}$.

The restriction property  given in Proposition \ref{P:Res} means that this does not 
occur for an observer using $P_E$.

\subsection{Anthropic estimates of the cosmological constant} \label{ss:BES}

This final example is based on  \cite{BES18, SPL}, which look  at a  multiverse, containing
a random number of universes, each with a random cosmological constant 
$\Lam \in [a,b]$. The proportion of mass in a universe which
turns into stars varies according to $\Lam$, and assuming that the number of observers 
is proportional to the number of stars  these papers ask if the value of $\Lam$ 
which we observe is likely.

We simplify the model in \cite{BES18}, in particular by restricting to a finite
number of AO. Let $n \ge 1$, $\Th=\{i/n, 0 \le i \le n-1\}$, $\kappa \in (0, n)$ and $p_n = \kappa/n$.
Let $m(t)$ be a continuous function on $[0,1]$, and let
$M_n = \sum_{\th \in \Th} n^{-1} m(\th)$. We assume that $M_n>0$ for all $n$.
We have $M_n \to M$ where
$$ M = \int_0^1 m(t) dt. $$
Let $(U^O_\th, \th \in \Th)$ be independent $\Ber(p_n)$ r.v. (For a justification 
of this `flat prior' in the cosmological context see \cite[p. 3737]{BES18}.)
We interpret $U^O_\th=1$ as the event that there exists a universe with 
cosmological constant $\th$. Let $n_0 \in \bN$ and 
let $(Y^O_\th, \th \in \Th)$ be independent  r.v. with values in $\{0,1, \dots, n_0 \}$
such that 
$$ \bE^O (Y^O_\th)=m(\th). $$
If $U^O_\th=1$ then $Y^O_\th$ is the number of AO in universe $\th$. 
Let $K = \Th \times \{0, \dots, n_0\}$; the set of AO is 
$$ \sX^O= \{ (\th,k): U^O_\th=1, 1\le k \le Y^O_\th \} \subset K. $$
Let $(\Om_O, \sF_O, \bP_O)$ be the probability space carrying 
$(U_\th^O, Y_\th^O, \th \in \Th)$.
As in Section 2 we construct the space $(\Om, \sF, \FOb, \bP)$, and write
$\sX$, $U_\th$, $Y_\th$ for the extensions of  $\sX^O$, $U^O_\th$, $Y^O_\th$ to this space.
Define $Z: K \to \bR$ by setting $Z((\th,k))=\th$; the observation of an AO (sited at location
$S \in \sX$) is that it is in a universe with cosmological constant $Z_S$. We define
$L_\th$ by \eqref{e:Lydef}, so that
$L_\th  = U_\th Y_\th. $
Thus $\bE(L_\th) = p_n m(\th)$, and 
$$  \bE(|\sX|) = \sum_{\th \in \Th} p_n m(\th) = \kappa M_n. $$

Set
\be
 \pi^{(O)}_n(\th) = n^{-1}, \q \pi^{(E)}_n(\th) = P_E(Z_S = \th), \q \pi^{(L)}_n(\th)= P_L(Z_S = \th).
 \ee
Thus $\pi^{(O)}_n$ gives the probability that a (uniformly randomly selected) universe has constant $\th$,
while $\pi^{(E)}_n$ and $\pi^{(L)}_n$ give the distribution of $\th$ for AO which use $P_E$ or $P_L$.
(One is tempted to call $\pi^{(O)}_n$ the prior, and $\pi^{(E)}_n$ and $\pi^{(L)}_n$
the posterior distributions of $\th$, but this is not exactly the standard Bayesian setup.)

By Lemma \ref{L:ZS}
\be
  \pi^{(E)}_n(\th) = \frac{ \bE(L_\th) }{ \bE (|\sX|)} = \frac{ m(\th)}{n M_n}. 
\ee 
Thus $\pi^{(E)}_n$ is weighted according to $m(\th)$. 
Note that this expression does not depend on $\kappa$.
This agrees with the calculation in \cite{BES18} -- see equation (14).

An interesting point is that if $F \in \sigma(U_t, Y_t, t \neq \th)$
then as $1_F$ and $L_\th$ are independent by Lemma \ref{L:ZS}
$$ P_E(F \cap \{Z_S = \th\}) = \bP( F) P_E( Z_S = \th). $$
So the AO has a flat distribution for the constants  of every universe other than its own.

If instead the AO uses $P_L$ then by Lemma \ref{L:ZPL}, 
\be \label{e:ZLB}
   P_L( Z_S = \th) = \bE \Big( {L_\th}{X^{-1}} \Big| X > 0 \Big)
   = \bP(X>0)^{-1} \bE( L_\th X^{-1} 1_{(X>0)} ). 
\ee

We now look at \eqref{e:ZLB} when $\kappa$ is either small, or very large.
To simplify calculations we assume further that
 $\bP( Y_\th =0) =0$ for all $\th$, and will take $n$ to be large. 
 Since we obtain different answers, it follows that $P_L(Z_S = \th)$ does depend on $\kappa$.
 We write $X = |\sX|$, $X_0 = X - U_\th Y_\th$.
 
 If $\kappa$ is small then both $\bP(X>0)$ and $\bP(X_0>0)$ are close to $\kappa$.
 We have, $\bP$-a.s.,
 $L_\th X^{-1} 1_{(X>0)}  =  U_\th  Y_\th/ (X_0+Y_\th) $,  
 so that 
 $$   U_\th 1_{(X_0=0)}  \le \frac{L_\th}{X}  1_{(X>0)}  \le  U_\th . $$
It follows that $\pi^{(L)}_n(\th)   \simeq  \pi^{(O)}_n(\th)$,
so that for an AO using $P_L$ the cosmological constant in its own universe still has a flat density.

Now let $\al \in (0,1)$ and let 
$\kappa$ be given by $\kappa_n = n^{1-\al}$. 
Then $q_n = \bP( X=0) = (1-p_n)^n \to 0$ as $ n \to \infty$, and
$$ \pi^{(L)}_n(\th) = \frac{1}{1-q_n} \bE \Big( \frac{ U_\th Y_\th}{X_0+Y_\th} \Big)
= \frac{1} {n(1-q_n)} \bE \Big( \frac{\kappa_nY_\th}{X_0 + Y_\th} \Big). $$
If $Y \ge 0$ then $\bE(1/Y) \ge 1/\bE(Y)$ and thus as $Y_\th$ and $X_0$ are independent,
$$  \bE \Big( \frac{\kappa_n Y_\th}{X_0 + Y_\th} \Big) \ge \bE \Big( \frac{\kappa_n Y_\th}{X_0 + n_0 } \Big) 
\ge \frac{ \kappa_n m(\th) }{\bE (X_0) + n_0} 
= \frac{ m(\th) }{ M_n + \kappa_n^{-1}(n_0 - m(\th))},$$
and thus $\liminf_n n \pi_{L,n}(\th) \ge m(\th)/M$.

Set $\mu_0 = \bE(X_0) = \kappa_n M_n - m(\th)$.
We have
$$ \bE \Big( \frac{\kappa_n Y_\th}{X_0 + Y_\th} \Big) \le
 \bE \Big( \frac{\kappa_n Y_\th}{X_0 + 1 } \Big)  = m(\th)  \bE\Big(  \frac{\kappa_n}{X_0 + 1 } \Big). $$
Straightforward computations give that $\bE( X_0- \mu_0 )^4 \le c \kappa_n^2$.
So if $\eps>0$ and $F_\eps=\{ X_0 < (1-\eps)\mu_0 \}$, then a fourth moment
estimate gives $ \bP( F_\eps) \le  c' \eps^{-4} \kappa_n^{-2}.$ 
Then
$$  \bE \Big( \frac{\kappa_n}{1+X_0} \Big) 
=    \bE \Big( \frac{\kappa_n 1_{F_\eps} }{1+X_0} \Big) +  \bE \Big( \frac{ \kappa_n1_{F^c_\eps} }{1+X_0} \Big) 
\le \kappa_n \bP(F_\eps) + \frac{\kappa_n}{ (1-\eps) \bE(X_0) }, $$
and it follows that
$$ \limsup_n \bE \Big( \frac{\kappa_n }{1+X_0} \Big) \le \frac{1}{ (1- \eps) M }. $$
Thus we have proved the following Lemma, which shows that in the large universe limit
$P_E$ and $P_L$ lead to the same estimates on $\th$.

\begin{lemma} If $\kappa_n = n^{1-\al}$ with $\al \in (0,1)$ then
\be
  \lim_n n \pi^{(L)}_n(\th) =  \lim_n n \pi^{(E)}_n(\th) = m(\th) /M. 
\ee
\end{lemma}

\begin{remark} {\rm
Many cosmological models involve infinite universes, which 
may then rise to an infinite number of observers.  It is therefore usual in the physics literature to look at
observers per unit volume --  see for example Section 5 of \cite{BES18}.
} \end{remark}
 
\sm {\bf Acknowledgment. } 
I wish to thank Richard A. Johns for telling me about the related Doomsday argument, 
Don Page for  pointing me to the paper \cite{HS}, and Chris Burdzy, 
Jason Swanson and Zoran Vondra\v cek for helpful conversations on this topic.
This research is partly supported by the Pacific Institute for the Mathematical Sciences.

\end{document}